\documentclass{article}

\usepackage{amsmath}
\usepackage{amsfonts}
\usepackage{amssymb}
\usepackage{graphicx}
\usepackage{makeidx}

\newtheorem{lemma}{{Lemma}}
\newtheorem{proposition}{{Proposition}}

\newtheorem{remark}{{Remark}}

\def\UUU{\mathcal{U}}
\def\gam{\gamma}
\def\lam{\lambda}

\def\vphi{\varphi}
\def\RRR{\mathbb{R}}

\def\E{\mathrm{E}}
\def\P{\mathrm{P}}
\def\C{\mathbb{C}}
\def\R{\mathbb{R}}
\def\O{{\cal O}}

\def\AA{\mathcal{A}}
\def\Z{\mathrm{Z}}
\def\Var{\mathrm{Var}}

\newcommand{\Vud}{{\mathrm{V}^{(u,d)}}}
\def\I{\;\mathbb{I}}
\newcommand{\A}[1]{#1\index{#1}}
\def\suchthat{\,{\Big |}\,}
\def\mstep{\bar{\mathrm{g}}}
\def\dstep{g}
\def\Summa{Q}
\def\Key{\mathrm{K}}
\def\Interval{L}
\def\coefficient{\mathrm{k}}
\def\F{{\cal F}}
\def\iVud{ {\buildrel {\circ \kern 4ex} \over \Vud} }

\title{A stochastic approximation algorithm with multiplicative step size adaptation}

\author{
\parbox{5cm}{\center Alexander Plakhov \\  plakhov@mat.ua.pt}
\hfill
\parbox{5cm}{\center Pedro Cruz \\ jpedro@mat.ua.pt } \\
Department of Mathematics\\
University of Aveiro --- Portugal
}

\makeindex

\begin{document}

\date{} 

\maketitle

\begin{abstract}
An algorithm of searching a zero of an unknown function $\vphi :
\, \R \to \R$ is considered, $\, x_{t} = x_{t-1} - \gamma_{t-1}
y_t$,\, $t=1,\ 2,\ldots$, where $y_t = \varphi(x_{t-1}) + \xi_t$
is the value of $\vphi$  measured  at $x_{t-1}$ with some error,
$\xi_t$ is this error. The step sizes $\gam_t > 0$ are random
positive values and are calculated according to the rule: $\,
\gamma_t = \min\{u\, \gamma_{t-1},\, \mstep\}$ if $y_{t-1} y_t >
0$, and $\gamma_t = d\, \gamma_{t-1}$, otherwise. Here $0 < d < 1
< u$,\, $\mstep
> 0$. The function $\vphi$ may have one or more zeros; the random values
$\xi_t$ are independent and identically distributed, with zero
mean and finite variance. Under some additional assumptions on
$\vphi$, $\xi_t$, and $\mstep$, the conditions on $u$ and $d$
guaranteeing a.s. convergence of the sequence $\{ x_t \}$, as well
as the conditions on $u$,\, $d$ guaranteeing a.s. divergence, are
determined. In particular, if $\P (\xi_1 > 0) = \P (\xi_1 < 0) =
1/2$ and $\P (\xi_1 = x) = 0$ for any $x \in \R$, it is
established that for $ud < 1$, convergence takes place, and for
$ud > 1$, divergence. Due to the multiplicative rule of updating
of $\gam_t$, it is natural to expect that $\{ x_t \}$ converges
rapidly: like a geometric progression (if convergence takes
place), but the limit value may not coincide with, but instead,
approximates one of zeros of $\vphi$. By adjusting the parameters
$u$ and $d$, one can reach necessary precision of approximation;
higher precision is obtained at the expense of lower convergence
rate.
\end{abstract}



\section{Introduction}

Consider the problem of finding a zero of a function $\vphi : \RRR
\to \RRR$. If there are several zeros, it is required to find at
least one of them. It is supposed that the function can be
measured at any point, with some random error. The standard
algorithm of stochastic approximation consists in calculating
successive approximations of the required value, $x_0$,\, $x_1$,\,
$x_2, \ldots$, according to the rule
\begin{equation}\label{eqal1}
x_{t} = x_{t-1} - \gamma_{t-1} y_t, \quad t=1,\ 2,\ldots,
\end{equation}
where
\begin{equation}\label{eqal2}
y_t = \varphi(x_{t-1}) + \xi_t \hspace{26mm}
\end{equation}
is the value of $\vphi$  measured at $x_{t-1}$,\, $\xi_t$ is the
measurement error; \,$\gam_0$,\, $\gam_1$,\, $\gam_2, \ldots$ is
the sequence of step sizes of the algorithm. Usually it is assumed
that the step sizes are positive real numbers satisfying the
relations $\sum \gam_t = \infty$,\, $\sum \gam_t^2 < \infty$.
Then, under some additional assumptions on $\vphi$ and $\xi_t$,
the algorithm a.s. converges to a zero point of $\vphi$ (see,
e.g., \cite{b007,b001}). In practice, however, the convergence rate of
this algorithm may prove to be unsatisfactory, therefore, when
solving practical tasks, various modifications of the algorithm
are used. There are widely utilized heuristical algorithms using
random, rather than deterministic, step size, which is corrected
in the course of the algorithm, according to the current data
\cite{a005,a011,a100,a006}. In particular, there is used the idea that prescribes to
decrease the step size if the sequence of increments $x_{t} -
x_{t-1}$ changes the sign often enough, indicating that the
current value $x_t$ is close to the set of zeros of $\vphi$, and
hence, the measurement error $\xi_t$ of the function is big enough
with respect to the function itself $\vphi(x_{t-1})$.
Alternatively, one should increase the step size, or leave it
unchanged. So, Kesten in the theoretical work \cite{a008}
considered an algorithm using (\ref{eqal1}), (\ref{eqal2}), and
the rule of modification of $\gam_t$:
\begin{equation}\label{eqK}
\gam_t = \gam(s_t), \ \ \ \ \ \ s_t = \left\{
    \begin{array}{lll}
    s_{t-1} & \textrm{ if } & y_{t-1} y_t > 0\\
    s_{t-1} + 1 & \textrm{ if } & y_{t-1} y_t \le 0,
    \end{array}
    \right.\quad t=2,3,\ldots.
\end{equation}
where $s_0 = 0$,\, $s_1 = 1$;\, $\gam(0)$,\, $\gam(1)$,\,
$\gam(2), \ldots$ is a sequence of positive numbers satisfying the
relations $\sum \gam(m) = \infty$,\, $\sum \gam^2(m) < \infty$.
Thus, the step size cannot increase in the course of algorithm; it
can only decrease or remain unchanged. It is supposed that there
is a unique zero of $\vphi$. Kesten proved that $x_t$ a.s.
converges to this zero point. A multidimensional version of this
algorithm is considered in \cite{a003}.

There are also heuristical procedures (in particular, in
artificial neural networks), where at each moment $t$ the step
size is multiplied by a positive constant less than 1, if the
measurement data indicate that $x_t$ is close enough to the zero
set of $\vphi$, and by a constant more than 1, elsewhere
\cite{a067,a001,a100,a101}. This kind of rules ensure
sufficiently high convergence rate, however the step size
converges like a geometric progression, therefore $\sum \gam_t <
\infty$, which means that the limit of $\{ x_t \}$ need not be a
zero point of $\vphi$, but instead, the sequence may "get stuck"
on its way to the set of zeros of $\vphi$. Nevertheless, such a
procedure may be justified if it gives a value close enough to one
of the zeros of $\vphi$.

In the present paper, a stochastic approximation algorithm
utilizing this rule of step size modification is considered.
Namely, the rule (\ref{eqal1}), (\ref{eqal2}), jointly with the
following rule
\begin{eqnarray}
  \gamma_t = \left\{
    \begin{array}{lll}
    \min\{u\, \gamma_{t-1},\, \mstep\} & \textrm{ if } & y_{t-1} y_t > 0,\\
    d\, \gamma_{t-1}                   & \textrm{ if } & y_{t-1} y_t \le 0,
    \end{array}
    \right.\quad t=2,3,\ldots. \label{eqal3}
\end{eqnarray}
is used. Here $0 < d < 1 < u$,\, $0 < \gamma_0$,\, $\gamma_1 \le
\mstep$,\,  $\mstep$ is a positive constant. Let us point out the
main differences between (\ref{eqal3}) and Kesten's rule
(\ref{eqK}). First, according to (\ref{eqal3}), $\gam_t$ can both
decrease and increase. Second, in Kesten's algorithm one always
has $\sum \gam_t = \infty$. On the other hand, it looks likely
that in the case of convergence of the algorithm (\ref{eqal1}),
(\ref{eqal2}), (\ref{eqal3}), $\gam_t$ converges like a geometric
progression (this conjecture will be justified in the section 3),
therefore the limit of algorithm may not be a zero point of
$\vphi$.

Suppose that $\{ \xi_t \}$ is a sequence of i.i.d.r.v. with zero
mean, besides $\P (\xi_t > 0) = \P (\xi_t < 0)$. Under some
additional assumptions on $\vphi$, $\xi_t$, and $\mstep$, stated
below, the process defined by (\ref{eqal1}), (\ref{eqal2}),
(\ref{eqal3}) a.s. diverges if $ud > 1$, and converges if $ud <
1$, moreover the limit of $\{ x_t \}$ belongs to $\UUU (\frac{\ln
u}{-\ln d})$. Here $\UUU(\lam)$,\, $0 < \lam < 1$, is a monotone
decreasing family of sets of real numbers, besides every set
$\UUU(\lam)$ contains the set $\Z$ of zeros of $\vphi$, and
$\partial (\UUU(\lam), \Z) \to 0$ as $\lam \to 1^-$. (Here by
definition $\partial (A,B) = \sup_{x\in A} \inf_{y\in B} |x-y|$
for any two sets of real numbers $A$ and $B$.)
 This statement is a consequence of the main
theorem, which will be stated in section 2 and proved in section
3. Thus, by adjusting the parameters $u$ and $d$ (for example,
fixing $u$ and letting $d \to 1/u - 0$), one can reach necessary
precision of the algorithm; higher precision is obtained at the
expense of lower convergence rate.

\section{Definition of the algorithm and statement of the main result}

Consider the algorithm given by (\ref{eqal1}), (\ref{eqal2}),
(\ref{eqal3}). The rule (\ref{eqal3}) means that at each instant
$t$, step size is multiplied by $u$ or by $d$, if the result of
multiplication is less than $\mstep$; otherwise, step size is set
to be $\mstep$. Thus, the maximal possible value of step size
equals $\mstep$.

The rule (\ref{eqal3}) can be written in the form
  \begin{equation}\label{eqal4n}
  \begin{array}{l@{ = }l}
 \ln \tilde \gamma_t & \ln \gamma_{t-1}
+ \ln u \cdot \I(y_{t-1} y_t>0) + \ln d \cdot \I(y_{t-1} y_t \le 0), \\
  \ln \gamma_t        & \min\{ \ln \tilde \gamma_t, \ln \mstep\}.
  \end{array}
  \end{equation}

Let us take the following assumptions:

\begin{description}

\item [A1] Denote $\F_t$, $t = 0,1,2,\ldots$ the $\sigma$-algebra
generated by $x_i$, $\gamma_i$, and $\xi_i$, $0\le i \le t$; then
$\xi_{t+1}$ does not depend on $\F_t$.

\item [A2] The values $\xi_t$ are identically distributed, with
zero mean and finite variance: $\E \xi_t = 0$,\, $\Var \xi_t =: S
< +\infty$.

\item [A3] (a) There exists $\Interval >0$ such that for any
interval
$I \subset [-\Interval,\, \Interval]$,  $\P(\xi_1 \in I)>0$;\\
\hspace*{-2mm}(b) $\P(\xi_1 = 0) = 0$.

\item [A4] $\varphi \in \C^1(\R)$ and $\sup_x |\varphi'(x)| =: M <
\infty$.

\item [A5] $\mstep < 2/M$.

\item [A6] There exists $R>0$ such that
\begin{itemize}
\item[(a)] $x \varphi(x) > 0$ as $|x| \ge R$, and \item[(b)]
$\displaystyle \inf_{|x|\ge R} \varphi^2(x) > \frac{\mstep M
S}{2-\mstep M}$.
\end{itemize}

\end{description}

\begin{remark}
From \A{A4} and \A{A6}\,(a) it follows that the set $\Z$ is
non-empty and is contained in $(-R,\, R)$.
\end{remark}

\begin{remark}
Note that assumptions \A{A4}--\A{A6} guarantee convergence of the
deterministic counterpart of algorithm (\ref{eqal1}),
(\ref{eqal2}), (\ref{eqal3}) (that is, of the algorithm with
$\xi_t \equiv 0$). Moreover, under these conditions, any
deterministic algorithm $x_t = x_{t-1} - \gamma_{t-1}
\varphi(x_{t-1})$ converges, whatever the sequence $\{\gamma_t\}$
satisfying $\gamma_t \le \mstep$.
\end{remark}

Introduce the functions:
 \begin{equation} \label{ast1}
 k_+(z) := \lim_{\epsilon\to0^+} \sup\{ \P( (\varphi_1+\xi_1)(\varphi_2+\xi_2) >
 0),\
 |\varphi_1 - z| < \epsilon,\ |\varphi_2 - z| < \epsilon\},
 \end{equation}
 \begin{equation} \label{ast2}
 k_-(z) := \lim_{\epsilon\to0^+} \inf\{ \P( (\varphi_1+\xi_1)(\varphi_2+\xi_2) >
 0),\
 |\varphi_1 - z| < \epsilon,\ |\varphi_2 - z| < \epsilon\};
 \end{equation}
one has $k_+(z) \ge 1/2$,\, $0 \le k_{\pm}(z) \le 1$,\,
$\lim_{z\to\infty} k_{\pm}(z)=1$.

Further, define the sets of real numbers
  \begin{equation}\label{ast3}
  V_{\pm}^{(a)} := \{ x : k_\pm(\varphi(x)) < a \}, \quad
  V_{\pm}^{[a]} := \{ x : k_\pm(\varphi(x)) \le a \};
  \end{equation}
obviously, $V_+^{(a)} \subset V_-^{(a)}$,\, $V_\pm^{(a)} \subset
V_\pm^{[a]}$ for any $a$.

Note that $V_+^{(a)}$ is open. Indeed, let $x \in V_+^{(a)}$, then
there exists $\epsilon > 0$ such that
$$
\sup\{ \P( (\varphi_1+\xi_1)(\varphi_2+\xi_2) > 0),\
 |\varphi_1 - \varphi(x)| < \epsilon,\ |\varphi_2 - \varphi(x)| <
 \epsilon\} =: c < a.
$$
Then for $x'$ close enough to $x$ one has $|\varphi(x')-
\varphi(x)| < \varepsilon/2$, hence
$$
\sup\{ \P( (\varphi_1+\xi_1)(\varphi_2+\xi_2) > 0),\
 |\varphi_1 - \varphi(x')| < \epsilon/2,\ |\varphi_2 - \varphi(x')| <
 \epsilon/2 \} \le c < a.
$$
This implies that $k_+(\varphi(x')) < a$, hence $x' \in
V_+^{(a)}$.

Denote also
  \begin{equation}\label{eq5}
  \coefficient := \frac{\ln(1/d)}{\ln(u/d)}.
  \end{equation}

Denote by $\Z$ the set of zeros of $\varphi$, i.e., $\Z:=\{ x :
\varphi(x)=0 \}$. Suppose that $x\in V_+^{(\coefficient)}$,\,
$x_{t-2} \in (x-\epsilon,\, x+\epsilon) \subset
V_+^{(\coefficient)}$, and $\gamma_{t-2} < \epsilon$, where
$\epsilon$ is a small positive number. Then, with a probability
close to 1,\, $x_{t-1}$ also belongs to a small (possibly larger)
neighborhood of $x$ contained in $V_+^{(\coefficient)}$, and
taking into account (\ref{ast1}) and (\ref{ast3}), one gets
  \[
  \begin{array}{l}
  \P( y_{t-1} y_t > 0 \suchthat |x_{t-2} - x| < \epsilon, \gamma_{t-2} < \epsilon ) =\\
  =\P( (\varphi(x_{t-2}) + \xi_{t-1}) (\varphi(x_{t-1}) + \xi_{t}) > 0
       \suchthat |x_{t-2}-x|<\epsilon, \gamma_{t-2} < \epsilon) <
       \coefficient.
  \end{array}
  \]
Then, using (\ref{eqal4n}) and (\ref{eq5}), one obtains
  \[
  \begin{array}{l}
  \E[ \ln \gamma_t-\ln \gamma_{t-1} \suchthat |x_{t-2}-x| < \epsilon, \gamma_{t-2} < \epsilon ]
  \le \\
  \ln u \cdot  \P(y_{t-1} y_t > 0 \suchthat |x_{t-2}-x| < \epsilon, \gamma_{t-2} < \epsilon ) +
  \ln d \cdot  \P(y_{t-1} y_t \le 0 \suchthat |x_{t-2} - x| < \epsilon, \gamma_{t-2} < \epsilon ) \\
  < \ln u \cdot \coefficient + \ln d \cdot (1-\coefficient) = 0.
  \end{array}
  \]

Thus, in a sense, the set $V_+^{(\coefficient)}$ can be regarded
to be a \textit{domain of decrease of step size}: if several
consecutive values of $x_t$ belong to $V_+^{(\coefficient)}$ and
are close enough to each other, and if the first term of the
sequence of corresponding step sizes $\gamma_t$ is small enough,
then the sequence of their mean values $E \gamma_t$ decreases.

Now, suppose that $x \in \R \setminus V_-^{[\coefficient]}$,\,
$x_{t-2} \in (x-\epsilon,\, x+\epsilon) \subset \R \setminus
V_-^{[\coefficient]}$, and that $\gamma_{t-2} < \epsilon$.
Analogously, for $\epsilon$ small enough, one has
  \[
  \P( y_{t-1} y_t > 0 \suchthat |x_{t-2} - x| < \epsilon, \gamma_{t-2} < \epsilon ) > \coefficient,
  \]
and then, using again (\ref{eqal4n}) and (\ref{eq5}) and taking
into account that for $\epsilon < \mstep/u^2$,\, $\tilde\gamma_t =
\gamma_t$, one obtains
  \[
  \begin{array}{l}
  \E[ \ln \gamma_t - \ln \gamma_{t-1} \suchthat |x_{t-2} - x| < \epsilon, \gamma_{t-2} < \epsilon] = \\
  \ln u \cdot \P( y_{t-1} y_t > 0 \suchthat |x_{t-2} - x | < \epsilon, \gamma_{t-2} < \epsilon]) 
 + \ln d \cdot \P( y_{t-1} y_t \le 0 \suchthat |x_{t-2} - x| < \epsilon, \gamma_{t-2} < \epsilon ])  \\
  > \ln u \cdot \coefficient + \ln d \cdot (1-\coefficient) = 0.
  \end{array}
  \]
Thus, the set $\R \setminus V_-^{[\coefficient]}$ can be regarded
as a \textit{domain of increase of step size}: if several
consecutive values of $x_t$ belong to $\R\setminus
V_-^{[\coefficient]}$ and are close enough to each other, and if
the first of the corresponding values of $\gamma_t$ is small
enough, then the sequence of their mean values $E \gamma_t$
increases.

Note that if $\coefficient > k_+(0)$ then, by virtue of
(\ref{ast3}), $\Z \subset V_+^{(\coefficient)}$, that is, all the
zeros of $\varphi$ belong to the region of decrease of step size.
On the other hand, if $\coefficient < \inf_z k_-(z)$ then
$V_-^{[\coefficient]} = \emptyset$, which means that the region of
increase of step size coincides with $\R$.

It seems likely that in the first case the algorithm can converge,
and in the second one, cannot. This conjecture is confirmed by the
following theorem, which is the main result of the paper.
 \vspace{2mm}

\textbf{Theorem} \textit{ Let the assumptions \A{A1}--\A{A6} be
satisfied; consider the process $\{ x_t,\ \gam_t \}$ defined by
(\ref{eqal1}), (\ref{eqal2}), (\ref{eqal3}). Recall that
$\coefficient = \frac{\ln(1/d)}{\ln(u/d)}$.
Then\\
(a) If $\coefficient > k_+(0)$ then $\{x_t\}$ a.s. converges to a
point from $V_-^{[\coefficient]}$.\\
(b) If $\coefficient  < \inf_z k_-(z)$ then $\{ x_t \}$ a.s.
diverges. }
 \vspace{2mm}

Suppose that $\P (\xi_1 = x) = 0$ for any real $x$ and that $\P
(\xi_1 > 0) = \P (\xi_1 < 0)$. Then the function $k(\cdot) :=
k_+(\cdot)$ coincides with $k_-(\cdot)$, is continuous, and is
given by
$$
k(z) = \P ((z + \xi_1)(z + \xi_2) > 0);
$$
$z = 0$ is the unique minimum of $k(\cdot)$, and $k(0) = \inf_z
k(z) = 1/2$. After a simple algebra, one can rewrite the
hypotheses of theorem in the form (a) $ud < 1$, (b) $ud > 1$.
Denote $\UUU(\lam) := V^{[\frac{1}{1+\lam}]} = \{ x :\,
k(\vphi(x)) \le \frac{1}{1 + \lam} \}$; \, $\UUU(\lam)$,\, $1 <
\lam < 1$ is a monotone decreasing family of sets containing $\Z$
and tending to $\Z$ as $\lam \to 1^-$.

Thus, one comes to  \vspace{2mm}

\textbf{Corollary} \textit{ Let, in addition to assumptions
\A{A1}--\A{A6}, $\P (\xi_1 = x) = 0$ for any $x \in \R$, and $\P
(\xi_1 > 0) = \P (\xi_1 < 0) = 1/2$. Consider the process defined
by (\ref{eqal1}), (\ref{eqal2}), (\ref{eqal3}). Then there exists
a monotone decreasing family of sets $\UUU(\lam)$,\, $0 < \lam <
1$ such that $\UUU(\lam) \supset
\Z$,\, $\partial (\UUU(\lam), \Z) \to 0$ as $\lam \to 1^-$, and\\
(a) if $ud < 1$ then $\{x_t\}$ a.s. converges to a point from
$\UUU (\frac{\ln u}{-\ln d})$;\\
(b) if $ud > 1$ then $\{x_t\}$ a.s. diverges. }

\begin{remark}
Theorem does not give any information about behavior of the
algorithm for the values $u$,\, $d$ such that
$$
\inf\nolimits_z k_-(z) \le \frac{\ln(1/d)}{\ln(u/d)} \le k_+(0).
$$
In particular, under the hypotheses of corollary, the case $ud =
1$ remains unexplored. These issues will be addressed elsewhere.
\end{remark}

\section{Proof of theorem}

First we prove 10 auxiliary lemmas, and then, basing on them, we
prove theorem.

Here all statements about random variables are supposed
to be true almost surely.

In the sequel, we shall mainly designate random values by Greek
letters, and real numbers and functions from $\R$ to $\R$, by
Latin ones; the letters $t$,\, $i$,\, $j$,\, $s$ will denote
integer non-negative numbers. The function $\vphi$ and the random
values $x_t$,\, $y_t$ are exceptions; also, traditional notation
$\epsilon$,\, $\delta$ for small positive numbers will be used.

\begin{lemma} 
If $\sum_t \gamma_t < \infty$ then the sequence $\{ x_t \}$ converges.
\end{lemma}

\textit{Proof.} Note that without loss of generality one can
assume that $x_0$ is bounded. Indeed, replacing $x_0$ by $\tilde
x_0 = x_0 \cdot \I(|x_0|<X)$ changes the process only with
probability $\P(|x_0|>X)$. By taking $X$ large enough, one can
make this probability arbitrarily small.

Let $C>0$; define the stopping time
$\tau_C = \inf \{ t : \sum_{i=0}^t \gamma_i > C\}$ and introduce the new
process $x_t^C$, $\gamma_t^C$ by
  \[
  \begin{array}{l}
  x_t^C = x_t,        \quad \gamma_t^C=\gamma_t \textrm{ as } t < \tau_c, \textrm{ and  } \\
  x_t^C = x_{\tau_C}, \quad \gamma_t^C=0 \textrm{ as } t \ge
  \tau_c.
  \end{array}
  \]
First, let us prove that the sequence $\{x_t^C\}$ is bounded.
Designate $M_R := \sup_{|x|\ge R} \frac{\varphi(x)}{x}$; from
\A{A4} it follows that $M_R<\infty$. One has
  \begin{equation}\label{eq9}
  |x_t^C| \le |x_{t-1}^C - \gamma_{t-1}^C \varphi(x^C_{t-1})| + \gamma_{t-1}^C |\xi_t|.
  \end{equation}
Using that $\gamma_{t-1}^C \le C$ and $|\varphi(x_{t-1})^C| \le |\varphi(0)| + M |x_{t-1}^C|$,
one obtains
  \begin{equation}\label{eq10}
  |x_t^C| \le |x_{t-1}^C|(1 + CM) + \gamma_{t-1}^C(|\varphi(0)|+|\xi_t|).
  \end{equation}

If $\gamma_{t-1}^C \le 2/M_R$, an even more precise estimate for
$x_t^C$ can be obtained. We shall distinguish between two cases:
(i) $|x_{t-1}|\le R$ and (ii) $|x_{t-1}^C| > R$.

In case (i), designating $\bar{b} := \sup_{|x|\le R}
|\varphi(x)|$, one has
  \begin{equation}\label{eq11}
  |x_{t-1}^C - \gamma_{t-1}^C \varphi(x_{t-1}^C)| \le |x_{t-1}^C| + \gamma_{t-1}^C \bar{b}.
  \end{equation}

In the case (ii) one has
  \[
  0 \le \gamma_{t-1}^C \frac{\varphi(x_{t-1}^C)}{x_{t-1}^C} \le \frac{2}{M_R} M_R = 2,
  \]
hence
  \begin{equation}\label{eq12}
  |x_{t-1}^C - \gamma_{t-1}^C\varphi(x_{t-1}^C)| \le |x_{t-1}^C|.
  \end{equation}
Thus, in both cases (i) and (ii), from (\ref{eq9}), (\ref{eq11}), and (\ref{eq12}) one gets
  \begin{equation}\label{eq13}
  |x_t^C| \le |x_{t-1}^C| + \gamma_{t-1}^C ( \bar{b} + |\xi_t| ).
  \end{equation}
The overall number of values of $t$ such that $\gamma_{t-1}^C \le 2/M_R$
is less than $CM_R/2$; therefore, using (\ref{eq10}) and (\ref{eq13}), one concludes that
  \begin{equation}\label{eq14}
|x_t^C| \le \left( |x_0| + \sum_{i=1}^t
\gamma_{i-1}^C(\bar{b}+|\varphi(0)|+|\xi_i|) \right) \cdot
(1+CM)^{CM_R/2}.
  \end{equation}

Denote $c_0 := \bar{b} + |\varphi(0)|+\E|\xi_1|$ and $\zeta_t :=
|\xi_t| - \E|\xi_t|$; using that $\sum_1^\infty \gamma_{i-1}^C \le
C$ one gets
  \begin{equation}\label{eq15}
|x_t^C| \le \left( |x_0| + C\, c_0 + \sum_{i=1}^t \gamma_{i-1}^C
\zeta_i \right)\cdot(1+CM)^{CM_R/2}.
  \end{equation}
Using that $\sum_1^\infty \E(\gamma_{t-1}^C \zeta_t)^2 = \E
\zeta_1^2 \cdot \sum_1^\infty \E(\gamma_{t-1}^C)^2 < \infty$, one
obtains that the martingale $\sum_1^t \gamma_{i-1}^C \zeta_i$ is
bounded; the value $x_0$ is also bounded, so, by (\ref{eq15}), one
concludes that the sequence $\{x_t^C\}$ is bounded.

Now, let us show that $\{x_t^C\}$ converges. From the definition of $x_t^C$ and
$\gamma_t^C$ it follows that
  \[
  x_t^C = x_0 - \sum_1^t \gamma_{i-1}^C \varphi(x_{i-1}^C) - \sum_1^t \gamma_{i-1}^C \xi_i.
  \]
Using that the sequence $\{ \varphi(x_{i-1}^C)\}$ is bounded and that $\sum_1^\infty \gamma_{i-1}^C \le C$,
one gets that the series $\sum_1^\infty \gamma_{i-1}^C \varphi(x_{i-1}^C)$ converges. Further, one has
  \[
\sum_1^\infty \E(\gamma_{t-1}^C \xi_t)^2 = S\cdot \sum_1^\infty
\E(\gamma_{t-1}^C)^2 < \infty,
  \]
hence the martingale $\sum_1^t \gamma_{i-1}^C \xi_i$ converges. This implies
that $\{x_t^C\}$ also converges.

Define the events $A_C = \{ \sum_t \gamma_t \le C\}$ and $A_\infty
= \{ \sum_t \gamma_t < \infty \}$. One has $A_\infty = \cup_C
A_C$. If $\sum_t \gamma_t \le C$ then $x_t^C=x_t$ for any $t$;
this means that $\I(A_C)\cdot(x_t^C - x_t)=0$ for any $t$ and $C$.
The sequence $\{ \I(A_C) x_t^C\}$ converges, therefore the
sequence $\{\I(A_C)x_t\}$ also converges, and passing to the limit
$C\to\infty$ one obtains that $\{ \I(A_\infty) x_t\}$ converges.
This means exactly that if $\sum_t \gamma_t < \infty$ then $\{ x_t
\}$ converges. \hfill$\Box$


\begin{lemma} 
If $\lim_{t\to\infty} x_t = x$ then $x \in V_-^{[\coefficient]}$.
\end{lemma}
\textit{Proof.} Note that, using \A{A3}\,(a), it is easy to show
that there exists $\delta_0 >0$ such that $\P(\xi_1 \not \in
[x-\Interval/2,\, x+\Interval/2]) > \delta_0$, whatever $x \in
\R$.

Next, for any $x\not \in V_-^{([\coefficient])}$ there exist
$w(x)>0$ and $0 < \epsilon(x) < \Interval/4$ such that the
following holds: for any two random variables $\phi_1$ and
$\phi_2$ satisfying the relations $|\phi_l - \varphi(x)| \le
\epsilon(x)$,\, $l=1,2$ one has
  \[
  \P( (\phi_1+\xi_1)(\phi_2+\xi_2) >0 ) > \frac{\ln(1/d)+w(x)}{\ln u + \ln(1/d)}.
  \]

Choose a countable set of intervals $U_i =
(\varphi(x_i)-\epsilon(x_i),\ \varphi(x_i)+\epsilon(x_i))$
covering the set $\varphi(\R\setminus V_-^{[\coefficient]})$, and
denote $w_i := w(x_i)$. Fix $i$ and $s\in\{0,\, 1,\, 2,\ldots\}$,
and define the auxiliary process $x_t^{(is)}$, $\gamma_t^{(is)}$
by formulas:
 \vspace{1mm}

if $t<s$ then $x_t^{(is)} = x_t$, \ and \ if $t \ge s$ then
  \begin{eqnarray}\label{eq16}
   x_t^{(is)} = \left\{
  \begin{array}{ll}
   x_{t-1}^{(is)} - \gamma_{t-1}^{(is)}\, y_t^{(is)} &
    \textrm{ if } \ \varphi(x_{t-1}^{(is)} - \gamma_{t-1}^{(is)}\, y_t^{(is)}) \in U_i,\\
   x_i & \textrm{ elsewhere};
   \end{array}
   \right.
 \end{eqnarray}
 \begin{equation}\label{eq17}
  y_t^{(is)} = \varphi(x_{t-1}^{(is)}) + \xi_t, \ \ \ \ \ \ \ \ \
  \ \ \ \ \ \ \ \ \ \ \  \ \ \hspace{35mm}
   \end{equation}
\begin{eqnarray}\label{eq18}
  \gamma_t^{(is)} = \left\{
  \begin{array}{l@{\textrm{ if }}l}
  \min\{u\gamma_{t-1}^{(is)},\, \mstep\} \ \ \ & \ y_{t-1}^{(is)}\, y_t^{(is)} > 0,\\
  d \gamma_{t-1}^{(is)}                \ \ \ & \ y_{t-1}^{(is)}\, y_t^{(is)} \le 0.
  \end{array}
  \right. \ \ \ \ \ \ \ \ \ \ \ \ \ \  \
  \end{eqnarray}

So, as $t \ge s$,\, $\varphi(x_t^{(is)})$ is forced to be
contained in $U_i$.

For $t\ge s+2$, using that $y_{t-1}^{(is)} =
\varphi(x_{t-2}^{(is)})+\xi_{t-1}$,\, $y_t^{(is)} =
\varphi(x_{t-1}^{(is)}) + \xi_t$,\, $\varphi(x_{t-2}^{(is)}) \in
U_i$, one obtains that
  \[
  \P( y_{t-1}^{(is)}\, y_t^{(is)} > 0) > \frac{\ln(1/d) + w_i}{\ln u + \ln(1/d)}
  \]
and
  \[
  \P( y_{t-1}^{(is)}\, y_t^{(is)} \le 0 ) < \frac{\ln u - w_i}{\ln u + \ln (1/d)},
  \]
hence
  $$
\displaystyle  \E[\ln u\cdot \I(y_{t-1}^{(is)}\, y_t^{(is)}>0)\,
+\, \ln d\cdot \I(y_{t-1}^{(is)}\, y_t^{(is)}\le 0) ] >
$$
$$
\displaystyle > \ln u \cdot \frac{\ln(1/d) + w_i}{\ln u +
\ln(1/d)}\ +\ \ln d \cdot \frac{\ln u - w_i}{\ln u + \ln(1/d)} =
w_i.
  $$

Consider variables $\phi_1=f_1(\xi_1, \xi_2)$ and $\phi_2 =
f_2(\xi_1,\xi_2)$ providing a solution of the (deterministic)
minimization problem:
  \[
  (\phi_1 + \xi_1)(\phi_2+\xi_2) \to \min,
  \]
subject to
  \[
  \begin{array}{l}
  |\phi_1 - \varphi(x_i)| \le \epsilon(x_i) \\
  |\phi_2 - \varphi(x_i)| \le \epsilon(x_i),\\
  \end{array}
  \]
and denote $Y_{t-1}^1 = f_1(\xi_{t-1},\xi_t) + \xi_{t-1}$,
$Y_{t}^2 = f_2(\xi_{t-1},\xi_t) + \xi_{t}$, $\eta_t = \ln u \cdot
\I(Y_{t-1}^1 Y_{t-1}^2 > 0) + \ln d \cdot \I( Y_{t-1}^1 Y_{t-1}^2
\le 0)$.
One has

(i) $\eta_t \le \ln u \cdot \I( y_{t-1}^{(is)}\, y_t^{(is)} > 0)
    + \ln d \cdot \I( y_{t-1}^{(is)}\, y_t^{(is)} \le 0)$;

(ii) $\eta_t$ are identically distributed, and $\E \eta_t \ge
w_i$;

(iii) the set of random variables $\{\eta_t,\ t \textrm{ even},\ t
\ge s+2 \}$ as well as the set $\{\eta_t,\ t \textrm{ odd},\ t \ge
s+2 \}$, are mutually independent.

From (ii)--(iii) it follows that almost surely $\sum_t \eta_t =
+\infty$, and from (i) it follows that
  \[
\sum_t [\ln u \cdot \I(y_{t-1}^{(is)}\, y_t^{(is)} > 0 ) + \ln d
\cdot \I( y_{t-1}^{(is)}\, y_t^{(is)} \le 0)]= +\infty,
  \]
so, by virtue of (\ref{eq18}), $\gamma^{(is)}$ does not go to zero.

Thus, there exists a random value $\chi > 0$ such that for
infinitely many values of $t$,\, $\gamma_t^{(is)} \ge \chi$.

Define a sequence of stopping times $\tau_0$, $\tau_1$, $\tau_2,
\ldots$ inductively, letting $\tau_0=0$ and $\tau_j = \inf\{ t >
\tau_{j-1} : \gamma_t^{(is)} \ge \chi\}$ for $j\ge 1$. The events
$B_j = \{ |\xi_{\tau_j+1} + \varphi(x_i)|>\Interval/2\}$ happen
with probability more that $\delta_0$ (recall the remark done in
the beginning of proof), and every event $B_j$, $j\ge 2$ does not
depend on the set of events $\{ B_1, \ldots, B_{j-1} \}$.
Therefore, for infinitely many values of $j$, $B_j$, takes place,
i.e., $|\xi_{\tau_j+1} + \varphi(x_i)| > \Interval/2$, and hence,
taking into account that $|y_{\tau_j+1}| \ge |\xi_{\tau_j+1} +
\varphi(x_i)| - |\varphi(x_{\tau_j}) - \varphi(x_i)|$ and
$|\varphi(x_{\tau_j})-\varphi(x_i)| < \epsilon(x_i) <
\Interval/4$, for these values of $j$ one has $|y_{\tau_j+1}| \ge
\Interval/4$. Thus, one concludes that
  \begin{equation}\label{eqast}
\textrm{for infinitely many values of }j, \ \ |\gamma_{\tau_j}
y_{\tau_j + 1}| \ge \chi\, \Interval/4.
  \end{equation}
Suppose that $x_t$ converges to a point from $\R\setminus
V_-^{[\coefficient]}$, then for some $i$ and $s$ one has $x_t \in
U_i$ as $t \ge s$, hence the process $x_t^{(is)}$,
$\gamma_t^{(is)}$ coincides with $x_t$,\, $\gamma_t$, and
therefore $\gamma_t\, y_{t+1} \to 0$ as $t \to \infty$. The last
relation contradicts (\ref{eqast}), thus Lemma 2 is proved.
\hfill$\Box$

\begin{lemma} 
Let $\sum_t \gamma_t = \infty$. Then for any open set $\O$
containing $\Z$ there exists a positive constant
$\dstep=\dstep(\O)$ such that either (i) for some $t$,\,
$x_t\in\O$, or (ii) for some $t$, $|x_t|<R$ and $\gamma_t >
\dstep$.
\end{lemma}
\textit{Proof.} Designate by $f$ the primitive of $\varphi$ such that $\inf_x f(x) = 0$.
Define the stopping time
  \[
  \tau = \tau(\O,\dstep) :=
  \inf \{ t : \textrm{ either
             (i) } x_t \in \O,
         \textrm{ or (ii) } |x_t|<R \textrm{ and } \gamma_t \ge \dstep\}.
  \]
The value of $\dstep \in (0,\mstep)$ will be specified below.

Consider the sequence $\E_t = \E[ f(x_t) \I(t<\tau)]$. Introducing
shorthand notation $f(x_t) =: f_t$,\, $\I(t<\tau)=: I_t$,\,
$f'(x_t)=:f_t'=\varphi_t$, and using that $I_t \le I_{t-1}$, one
gets
  \begin{equation}\label{eq20}
  E_t - E_{t-1} = \E[f_t \I_t - f_{t-1} \I_{t-1}]\, \le\, \E[(f_t - f_{t-1}) \I_{t-1}].
  \end{equation}
Next, we utilize the Taylor decomposition
$$
f_t = f(x_{t-1} - \gamma_{t-1} y_t) = f_{t-1} - f'_{t-1}\,
\gamma_{t-1} y_t + \frac 12\, f''(x')\, \gamma_{t-1}^2 y_t^2,
 $$
$x'$ being some point between $x_{t-1}$ and $x_t$. Substituting
$y_t = \varphi_{t-1}+\xi_t$ and recalling that
$f'_{t-1}=\varphi_{t-1}$ and $f''(x')=\varphi'(x') \le M $, one
obtains
  \begin{equation}\label{eq21}
f_t - f_{t-1} \le -\gamma_{t-1}\, \varphi_{t-1} (\varphi_{t-1} +
\xi_t) + {M\over 2}\, \gamma_{t-1}^2\, (\varphi_{t-1} + \xi_t)^2.
  \end{equation}
Using (\ref{eq20}) and (\ref{eq21}) and taking into account that
each of the values $\gamma_{t-1}$,\, $\varphi_{t-1}$,\, $\I_{t-1}$
is mutually independent with $\xi_t$ (see \A{A1}), one gets
  \begin{equation}\label{eq23}
  \begin{array}{l}
  E_t - E_{t-1} \le \E[
           (-\gamma_{t-1}\, \varphi_{t-1}^2 -
       \gamma_{t-1}\, \varphi_{t-1} \xi_t +
       {M\over 2} \gamma_{t-1}^2\, \varphi_{t-1}^2 +
       M\gamma_{t-1}^2\, \varphi_{t-1} \xi_t +
       {M\over 2} \gamma_{t-1}^2\, \xi_t^2) \I_{t-1}] =\\
        = \E[ (-\varphi_{t-1}^2 + \frac M2 \gamma_{t-1}\, \varphi_{t-1}^2 + {M\over 2} \gamma_{t-1} S) \gamma_{t-1} \I_{t-1}] = \\
    = \E[ (-\varphi_{t-1}^2(1-M\gamma_{t-1}/2) + M\gamma_{t-1} S/2) \gamma_{t-1} \I_{t-1}].
  \end{array}
  \end{equation}

If $\I_{t-1}=1$ then~ either~ (i)~ $x_{t-1} \in [-R,R]\setminus\O$
and $\gamma_{t-1} < \dstep$,~~ or~ (ii)~ $|x_{t-1}| \ge R$.

In the case (i) one has
\begin{equation}\label{eq23.1}
  -\varphi_{t-1}^2(1-M\gamma_{t-1}/2) + M\gamma_{t-1}S/2 \le -c_0(1-M\dstep/2) + M \dstep S/2 =: -c'_\dstep,
 \end{equation}
where $c_0 := \inf\{ |\varphi(x)| : x\in[-R,R]\setminus\O \}$;
obviously, $c_0>0$. Let us fix a $\dstep \in (0,\mstep)$ such that
$c'_\dstep>0$.

In the case (ii), designating $b_0 := \inf_{|x|\ge R}
\varphi^2(x)$, one has
\begin{equation}\label{eq23.2}
  -\varphi_{t-1}^2(1-M\gamma_{t-1}/2) + M\gamma_{t-1} S/2 \le
   -b_0(1-M\mstep/2)+M\mstep S/2 =: -c''.
 \end{equation}
Using \A{A6}, one gets that $c''>0$.

Denote $c=\min\{c'_\dstep,c''\}$. The relations (\ref{eq23.1}) and
(\ref{eq23.2}) imply that if $\I_{t-1}=1$ then
$-\varphi_{t-1}^2(1-M\gamma_{t-1}/2) + M\gamma_{t-1} S/2 \le -c <
0$, hence, by virtue of (\ref{eq23}),
  \begin{equation}\label{eq24}
  E_t - E_{t-1} \le -c \cdot \E[\gamma_{t-1} \I_{t-1}].
  \end{equation}

Summing up both sides of (\ref{eq24}) over $t=1,\ldots,s$ and
denoting $\I_\infty = \I(\tau=\infty)=\min_t \I_t$, one obtains
  \[
  \E_s - \E_0 \le -c \cdot \E \left[\sum_{i=0}^{s-1} \gamma_i \cdot \I_\infty
  \right].
  \]
One has $\E_s \ge 0$, and $x_0$ is bounded, hence $E_0 < \infty$.
Thus, for arbitrary $s$
  \[
  \E \left[\sum_{i=0}^{s-1} \gamma_i \cdot \I_\infty \right] \le \frac{\E_0}{c} < \infty.
  \]
This implies that a.s. either $\sum_0^\infty \gamma_i < \infty$,
or $\tau = \infty$. Lemma 3 is proved. \hfill$\Box$
 \vspace{2mm}

Denote $c_1 := 1- M \mstep/2$. Recall that $f$ is the primitive of
$\varphi$ such that $\inf_x f(x) = 0$; the assumption \A{A6}
implies that $\lim_{x\to\pm \infty} f(x) = +\infty$. Denote $H:=
\sup_{|x|\le R} f(x)$. Denote also~ $c_3:= \mstep \cdot \sup\{
|\varphi(x)| : f(x) \le H\} + 1$,~ $z^{l} := \inf \{ x : f(x) \le
H \} - c_3$,\ $z^{r} := \sup\{ x : f(x) \le H\} + c_3$,~ $c_2 :=
\inf\{ |\varphi(x)| : x \in [z^{l},\, z^{r}] \setminus\O\}$,~ and~
$\Key := \sup\{ |\varphi(x)| : x\in[z^{l},\, z^{r}]\}$.~
Obviously, $c_1
> 0$ and $\Key\ge c_2 > 0$.

Fix an open set $\O$ containing $\Z$. Let $\dstep > 0$,\, $0 < w <
1$. We shall say that a (finite or infinite) deterministic
sequence $\{z_0, z_1, z_2, \ldots \}$ is $(\dstep,\,
w)$-admissible if $|z_0|\le R$ and there exist deterministic
sequences $\{q_t\},$ $\{h_t\}$ such that

1)   $|h_t| \le w$;

2)  if $\{ z_0,z_1,\ldots,z_t\}
\subset[z^{l},\,z^{r}]\setminus\O$~ then~     $\dstep d^2 \le q_s
\le \mstep$,~ $s=0,1,\ldots,t$;

3) $z_t = z_{t-1} - q_{t-1}\, \varphi(z_{t-1})-h_t$,~
$t=1,2,\ldots$.


\begin{proposition} 
There exists constants $t_0$ and $w$ such that any $(\dstep,\,
w)$-admissible sequence $\{z_t,\ t=0,\, 1,\ldots,t_0\}$ has
non-empty intersection with $\O$.
\end{proposition}
\textit{Proof.} Let $w:= \min\{ 1,\, \dstep d^2 c_2^2
c_1/(2\Key)\}$. Designate $\tilde{t}=\inf\{ t:z_t\in\O\}$;\,
$\tilde{t}$ takes values from $\{0,\, 1, \ldots, t_0,\,
+\infty\}$. We shall use shorthand notation $f_t := f(z_t)$,
$f'_t= \varphi_t := \varphi(z_t)$. One has
  \begin{equation}\label{eq25}
f_t = f(z_{t-1}-q_{t-1} \varphi_{t-1} - h_t) = f(z_{t-1} - q_{t-1}
\varphi_{t-1}) - f'(\tilde z).h_t,
  \end{equation}
where $\tilde z$ is a point between $z_{t-1} - q_{t-1}
\varphi_{t-1}$ and $z_{t-1} - q_{t-1} \varphi_{t-1} -h_t$.

Next, one has
  \begin{equation}\label{eq26}
f(z_{t-1}-q_{t-1} \varphi_{t-1}) = f_{t-1} - f'_{t-1} q_{t-1}
\varphi_{t-1} + \frac 12 f''(\hat z)\, q_{t-1}^2 \varphi_{t-1}^2,
  \end{equation}
where $\hat z$ is a point between $z_{t-1}$ and $z_{t-1}-q_{t-1}
\varphi_{t-1}$.

We are going to prove by induction that
  \begin{equation}\label{eq29}
  \textrm{if } 0\le s \le \tilde{t} \ \textrm{ then } \ f_s \le H-s\cdot \dstep d^2 c_2^2 c_1 / 2.
  \end{equation}
For $s=0$,~ (\ref{eq29}) follows from the condition $|z_0|\le R$
and the definition of $H$. Now, let $1 \le t \le \tilde{t}$;~
suppose that formula (\ref{eq29}) is true for $0 \le s \le t-1$
and prove it for $s = t$. For $0 \le s \le t-1$, one has $f(z_s)
\le H$,\, $z_s \not \in \O$, therefore $z_s \in [z^{l},\, z^{r}]
\setminus \O$; hence, by virtue of 2), $\dstep d^2 \le q_s \le
\mstep$ for $0 \le s \le t-1$.
 One has $f(z_{t-1}) \le H$,\, $|q_{t-1} \varphi_{t-1}| \le \mstep
\cdot \sup\{ |\varphi(x)| : f(x) \le H\}$, and $|h_t| \le w \le
1$, hence $|q_{t-1} \varphi_{t-1}| \le c_3$,\, $|q_{t-1}
\varphi_{t-1} + h_t| \le c_3$, and so, $z_{t-1} - q_{t-1}
\varphi_{t-1} \in [z^{l},\, z^{r}]$,\,
$z_{t-1}-q_{t-1}\varphi_{t-1} -h_t \in [z^{l},\, z^{r}]$, thus
$\tilde z$ also belongs to $[z^{l},\, z^{r}]$. This implies that
$|\varphi(\tilde z)| = |f'(\tilde z)| \le \Key$. Then, combining
(\ref{eq25}) and (\ref{eq26}) and using that $|h_t| \le w$ and
$|f''(\hat z)| = |\varphi'(\hat z)|\le M$, one obtains
  \begin{equation}\label{eq27}
  f_t \le f_{t-1} - q_{t-1} \varphi^2_{t-1}(1-{1\over 2} q_{t-1} M) + w\Key.
  \end{equation}
One has $z_{t-1} \in [z^{l},\, z^{r}] \setminus \O$, hence
$|\varphi(z_{t-1})| = |\varphi_{t-1}| \ge c_2$. Using also that
$q_{t-1} \ge \dstep d^2$,\, $1-{1 \over 2} q_{t-1} M \ge c_1$, and
$w\Key \le \dstep d^2 c_2^2 c_1/2$, one gets from (\ref{eq27})
that
  \[
  f_t \le f_{t-1} - \dstep d^2 c_2^2 c_1 / 2,
  \]
and using the induction hypothesis, one concludes that
  \[
  f_t \le H - t \cdot \dstep d^2 c_2^2 c_1 /2.
  \]

Formula (\ref{eq29}) is proved.

Let $t_0 := \lfloor 2H/(\dstep d^2 c_2^2 c_1) \rfloor  + 1$; here
$\lfloor z \rfloor$ stands for the integral part of $z$. Then,
taking into account that $f_s \ge 0$, from (\ref{eq29}) one
concludes that $\tilde{t} < t_0$, thus Proposition 1 is proved.
\hfill$\Box$.


\begin{proposition} 
If $\gamma_{t-1} < 1/(3M)$,\, $|\xi_t|<c_2$,\, $|\xi_{t+1}|<
c_2$,\, $x_{t-1}$ and $x_t$ belong to $[z^{l},\, z^{r}] \setminus
\O$,~ then $\gamma_{t+1} \ge \gamma_t$.
\end{proposition}
\textit{Proof.} Using notation $\varphi_t := \varphi(x_t)$, one gets
$$
\varphi_t = \varphi(x_{t-1} - \gamma_{t-1}(\varphi_{t-1}+\xi_t)) =
 \varphi_{t-1} - \varphi'(\tilde x) \cdot
 \gamma_{t-1}(\varphi_{t-1}+\xi_t),
$$
 where $\tilde x$ is a point between $x_{t-1}$ and $x_t$.
Therefore,
  \[
\varphi_{t-1} \varphi_t = \varphi^2_{t-1} \cdot[ 1-
\varphi'(\tilde x) \gamma_{t-1} \cdot(1+\xi_t / \varphi_{t-1})].
  \]
Using that $|\varphi'(\tilde x)| \le M$,\, $\gamma_{t-1} <
1/(3M)$,\, $|\xi_t| < c_2$,\, $|\varphi_{t-1}| \ge c_2$, one
obtains $1-\varphi'(\tilde x)\, \gamma_{t-1} \cdot (1 + \xi_t /
\varphi_{t-1}) \ge 1/3$, hence $\varphi_{t-1} \varphi_t >0$.
Further, using that $|\xi_t| < c_2$,\, $|\xi_{t+1}|<c_2$,\,
$|\varphi_{t-1}| \ge c_2$,\, $|\varphi_t| \ge c_2$, one gets
  \[
y_{t}\, y_{t+1} = \varphi_{t-1} \varphi_t
\cdot(1+\xi_t/\varphi_{t-1})(1+\xi_{t+1}/\varphi_t) > 0.
  \]
This implies that $\gamma_{t+1} = \min\{ u \gamma_t, \mstep\} \ge
\gamma_t$. \hfill$\Box$


\begin{lemma}
For any open set $\O$, containing $\Z$, and any $\dstep > 0$~
there exists $\delta = \delta(\O, \dstep) > 0$ such that
  \[
   \text{if } \ |x_0| \le R, \ \gamma_0 \ge \dstep \ \text{ then
 } \ \ \P(\textrm{for some } t, \ x_t \in \O) \ge \delta.
  \]
\end{lemma}
\textit{Proof.} Without loss of generality suppose that $\dstep <
1/(3M)$. Define the event
  \[
  A := \{ |\xi_i| < \min\{ c_2,\, w/\mstep \}, \ i=1,2,\ldots,t_0\},
  \]
where $w$ and $t_0$ are the same as in the proof of Proposition
1:~ $w = \min\{ 1,\, \dstep d^2 c_2^2 c_1/(2\Key)\}$,\ $t_0 =
\lfloor 2H/(\dstep d^2 c_2^2 c_1) \rfloor  + 1$.

Denote
  \[
  \delta := P(A) = ( \P( |\xi_1| < \min\{ c_2,\, w/\mstep\}))^{t_0};
  \]
by virtue of \A{A3}\,(a), $\delta > 0$. Let us show that for any
elementary event $\omega \in A$, the sequence $\{ z_t =
x_t(\omega),\ t=0, 1, \ldots, t_0\}$ is $(\dstep,\,
w)$-admissible.

One has $|z_0|=|x_0(\omega)| < R$. Further, one has $z_t = z_{t-1}
- q_{t-1} \varphi(z_{t-1})- h_t$, with $q_{t-1} =
\gamma_{t-1}(\omega)$,\, $h_t = \gamma_{t-1}(\omega)\,
\xi_t(\omega)$, and using that $\gamma_{t-1}(\omega) \le \mstep$
and $|\xi_t(\omega)| < \omega / \mstep$, one gets $|h_t| \le w$.
Thus, conditions 1) and 3) are verified.

Now, let $\{z_0, z_1, \ldots,z_t\} \subset [z^{l},\, z^{r}]
\setminus \O$,~ $t \le t_0$. Let $s_0 \in \{ 0,1,2,\ldots,t\}$ be
the minimal value such that $q_{s_0} = \min\{ q_0, q_1, \ldots,
q_t \}$. If $s_0=0$ then $\min\{ q_0, q_1,\ldots, q_t \} = q_0 =
\gamma_0(\omega)\ge \dstep \ge \dstep d^2$.
 If $s_0=1$ then
$\min\{ q_0, q_1, \ldots, q_t\} = q_1 = \gamma_1(\omega) \ge
\dstep d \ge \dstep d^2$. If $s_0 \ge 2$ then
$\gamma_{s_0-2}(\omega) \ge 1/(3M)$; otherwise, using that
$|\xi_{s_0-1}| < c_2$,\, $|\xi_{s_0}| < c_2$,\,
$x_{s_0-2}(\omega)$ and $x_{s_0-1}(\omega)$ belong to $[z^{l},\,
z^{r}] \setminus \O$, and applying Proposition 2, one would
conclude that $\gamma_{s_0}(\omega) \ge \gamma_{s_0-1} (\omega)$,
which contradicts the definition of $s_0$.

Thus, $\gamma_{s_0}(\omega) \ge 1/(3M) \cdot d^2 \ge \dstep d^2 $,
and therefore, $\min\{ q_0, q_1, \ldots, q_t \} =
\gamma_{s_0}(\omega) \ge \dstep d^2$. So, the condition 2) is also
verified.

Now, applying Proposition 1 to the $(\dstep,\, w)$-admissible
sequence $\{z_t\}$, one concludes that there exists a non-negative
$\tau \le t_0$ such that $z_{\tau} = x_\tau(\omega) \in \O$. This
implies that
  \[
  \P(\textrm{for some } t, \ x_t \in \O) \ge \P(A) = \delta.
  \]
\hfill$\Box$


\begin{lemma}
If $\sum_t \gamma_t = \infty$ then for any open set $\O$
containing $\Z$ there exists $t$ such that $x_t \in \O$.
\end{lemma}
\textit{Proof.} Let us fix an open set $\O \supset \Z$, and denote
$\delta = \delta(\O, \dstep(\O))$. Combining Lemma 3 and Lemma 4,
one concludes that for any $\O \supset \Z$ there exists $\delta >
0$ such that whatever the initial conditions $x_0$, $\gamma_0$,
$\gamma_1$,
  \[
\P(\textrm{for some } t, \ x_t \in \O \suchthat \sum_t \gamma_t =
\infty) > \delta.
  \]
Then one can choose a measurable integer-valued function
$n(\cdot,\cdot,\cdot)$ defined on $\R \times (0,\mstep] \times
(0,\mstep]$ such that for $\nu=n(x_0,\gamma_0,\gamma_1)$ one will
have
  \[
  \P(\textrm{for some } t\le \nu, \ x_t \in \O \suchthat \sum_t \gamma_t=\infty) > \delta/2
  \]

Designate
  \[
\bar p = \sup \P(\textrm{for all } t, \ x_t \not \in \O \suchthat
\sum_t \gamma_t = \infty),
  \]
the supremum being taken over all the initial conditions $x_0$, $\gamma_0$, $\gamma_1$.
Fix $x_0$, $\gamma_0$, $\gamma_1$, then
  \begin{equation}\label{eq30}
  \begin{array}{l}
  \P(\textrm{for all } t, \ x_t \not \in \O \suchthat \sum_t \gamma_t = \infty ) =\\
  = \P(\textrm{for all } t > \nu, \ x_t \not \in \O \suchthat
\textrm{for all } t \le \nu, \ x_t \not \in \O \textrm{ and } \sum_t \gamma_t =\infty) \cdot \\
\cdot \P(\textrm{for all } t \le \nu, \ x_t \not \in \O\, | \sum_t
\gamma_t = \infty) \le \bar p\, (1-\delta/2).
  \end{array}
  \end{equation}

Taking supremum of the left hand side of (\ref{eq30})  over all
$(x_0, \gamma_0, \gamma_1) \in \R \times (0,\mstep] \times
(0,\mstep]$, one obtains $\bar p \le \bar p\, (1- \delta/2)$,
hence $\bar p = 0$. Lemma~5 is proved. \hfill$\Box$.
 \vspace{2mm}

Denote $\O_* = \{ x: |\varphi(x)| < \Interval/2 \}$.

\begin{lemma}
 For any open bounded sets $\mathcal{O}$,\, $\mathcal{O}_1$
such that $\bar{\mathcal{O}} \subset \mathcal{O}_1 \subset \O_*$
and for any $w > 0$ there exists $\delta = \delta(\O, \O_1, w)> 0$
such that
 $$
\text{if } \ x_0 \in \mathcal{O} \text{ then } \ \P(\textrm{for
some } n,\ x_{n} \in \O_1 \text{ and } \gamma_{n} < w) \ge \delta.
 $$
\end{lemma}

\textit{Proof.} Denote $n = \lfloor \frac{\ln\mstep - \ln w}{\ln
(1/d)} \rfloor + 2$. Denote also
 $$
\varepsilon = \min \left\{ \frac{\Interval}{2},\ \frac{
\partial (\mathcal{O},\, \mathbb{R} \setminus \mathcal{O}_1)}{n \mstep}
\right\},
 $$
where $\partial (A,B) := \sup_{x\in A} \inf_{y\in B} |x-y|$ for
arbitrary sets of real numbers $A$,\, $B$. Using assumption
\A{A3}\,(a), one obtains that there exists $\delta_1 > 0$ such
that for any $x \in \O_1$ and for any integer $t$,
 $$
\P \left( (-1)^{t-1} \varphi(x) < (-1)^t \xi_1 < (-1)^{t-1}
\varphi(x) + \varepsilon \right) \ge \delta_1.
 $$
This implies that if $x_0 \in \O$ then
 $$
\P (0 < (-1)^t y_t < \varepsilon,\ \text{dist}(x_{t-1},\, \O) < (t
- 1) \mstep \varepsilon,\ t = 1,\, 2,\ldots, n+1) \ge
\delta_1^{n+1}.
 $$
Denoting $\delta = \delta_1^{n+1}$, one concludes that the
following statements (i) and (ii) hold with probability at least
$\delta$:

(i) dist$(x_n,\, \O) < n \mstep \varepsilon \le$ dist$(\O,\,
\mathbb{R} \setminus \O_1)$, hence $x_n \in \O_1$;

(ii) as $t = 2,\, 3,\ldots, n+1$, one has $y_{t-1} y_t < 0$, hence
$\gamma_t = d \gamma_{t-1}$, therefore

$\gamma_{n} = d^{n-1} \gamma_1 \le d^{n-1} \mstep < w$.\\
Lemma~6 is proved. \hfill$\Box$

\begin{lemma}
If $\sum_t \gamma_t = \infty$,\ $\O$ is an open set containing
$\Z$, and $w > 0$ then for some $t$,~ $x_{t-1} \in \O$ and
$\gamma_t < w$.
\end{lemma}

\textit{Proof.} Without loss of generality, suppose that $\O$ is
bounded and $\O \subset \O_*$. Choose an open set $\O_1$ such that
$\Z \subset \O_1$,\, $\bar\O_1 \subset \O$;~ applying Lemmas 5 and
6, one gets that for $\delta = \delta(\O_1, \O, w)$ and for
arbitrary initial conditions,
 $$
\P (\textrm{for some } t,\ x_{t} \in \O \textrm{ and } \gamma_t <
w)
> \delta.
 $$
Repeating the argument of Lemma 5, one concludes that  there
exists $t$ such that $x_{t} \in \O$ and $\gamma_t < w$.
     \hfill$\Box$
 \vspace{2mm}

From now on we suppose that $\coefficient > k_+(0)$. Choose $k'$
such that $k_+(0) < k' < \coefficient$; using \A{A3}\,(b), one
obtains that for some $\varepsilon_0 > 0$,\, $\P ( \xi_1 \xi_2 >
0, \text{ or } |\xi_1| < \varepsilon_0, \text{ or } |\xi_2| <
\varepsilon_0) \le k'$. Denote $\O_0 = \{ x:\, |\varphi(x)| <
\varepsilon_0 \}$ and $\tau = \inf \{ t: \ x_t \not\in \O_0 \}$.
Without loss of generality, suppose that $\O_0$ is bounded.

\begin{lemma}
Suppose that $\coefficient > k_+(0)$,~ then there exist a constant
$b
> 0$ and a monotone decreasing function $p(\cdot)$ such that
$\lim_{a\to+\infty} p(a) = 0$ and
 $$
\text{if } \ \gamma_0 < w \text{ then } \ \P (\ln\gamma_t < \ln v
- bt \text{ for all } t < \tau ) > 1 - p(v/w).
 $$
\end{lemma}

\textit{Proof.} Define the sequences $\{ \rho_t \}$ and $\{
\sigma_t \}$ by
\begin{eqnarray*}
\rho_t &=& \ln u \cdot \I(\xi_{t-1} \xi_t > 0, \text{ or }
|\xi_{t-1}|
< \varepsilon_0, \text{ or } |\xi_t| < \varepsilon_0) +\\
 &+& \ln d \cdot \I(\xi_{t-1} \xi_t \le 0 \ \, \& \, \ |\xi_{t-1}| \ge
\varepsilon_0 \ \, \& \ \, |\xi_t| \ge \varepsilon_0),
\end{eqnarray*}
$$
\sigma_t = \ln w + \sum_{i=1}^t \rho_i.
$$
Using (\ref{eqal4n}) and definition of $\tau$, one obtains that
for all $t < \tau$,\, $\gamma_t \le \sigma_t$. The variables
$\rho_t$ are identically distributed, take the values $\ln u$ and
$\ln d$, and
\begin{eqnarray*}
E\rho_t &=& \ln u \cdot \P(\xi_{t-1} \xi_t > 0, \text{ or }
|\xi_{t-1}| < \varepsilon_0, \text{ or } |\xi_t| < \varepsilon_0) +\\
 &+& \ln d \cdot \P(\xi_{t-1} \xi_t \le 0 \ \, \& \, \ |\xi_{t-1}| \ge
\varepsilon_0 \ \, \& \ \, |\xi_t| \ge \varepsilon_0) \le\\
&\le& \ln u \cdot k' + \ln d \cdot (1 - k') < \ln u \cdot
\coefficient + \ln d \cdot (1-\coefficient) = 0.
\end{eqnarray*}
Moreover, the variables in the set $\{ \rho_t, \ t \text{ even}
\}$, as well as the variables in the set $\{ \rho_t, \ t \text{
odd} \}$, are independent.

Denote $b = -E\rho_t/2$. One has
$$
\P (\ln\gamma_t < \ln v - b t \ \text{ for all } t < \tau ) \ge \P
(\sigma_t < \ln v - b t \ \text{ for all } t) =
$$
$$
= \P (\sum_{i=1}^t (\rho_i + 2b) < \ln v - \ln w + b t \ \text{
for all } t ) \ge 1 - p(v/w),
$$
where $p(a) = p_1(a) + p_2(a)$,
 $$
p_1(a) = \P \left( {\sum_{1\le i\le t}}' (\rho_i + 2 b) \ge  \frac
{\ln a}2 +  \frac b2\, t \ \text{ for all } t \right),
 $$
 $$
p_2(a) = \P \left( {\sum_{1\le i\le t}}'' (\rho_i + 2 b) \ge \frac
{\ln a}2 +  \frac b2\, t \ \text{ for all } t \right);
 $$
the sum $\sum'$ ($\sum''$) is taken over the even (odd) values of
$i$. Both $\sum'$ and $\sum''$ are sums of i.i.d.r.v. with zero
mean, hence both $p_1(a)$ and $p_2(a)$ tend to zero as $a \to
+\infty$. Lemma 8 is proved.
 \hfill$\Box$

Define the stopping times $\tau_v = \inf \{ t: \ x_t \not\in \O_0
\text{ or } \ln\gamma_t \ge \ln v - bt \}$. Recall that $f$ is the
primitive of $\varphi$ such that $\inf_x f(x) = 0$. Fix an open
set $\O'$ such that $\Z \subset \O' \subset \O_0$ and
$\sup_{x\in\O'} f(x) < \inf_{x\not\in\O_0} f(x)$, and denote
$\delta = \inf_{x\not\in\O_0} f(x) - \sup_{x\in\O'} f(x)$.

\begin{lemma}
Let $\coefficient > k_+(0)$,\ $x_0 \in \O'$, and $\gamma_0 < w$,
then
 $$
\P (\tau_v < \infty ) \le K\, v^2 + p(v/w);
 $$
here $K$ is a positive constant, and $p(\cdot)$ satisfies the
statement of lemma 8.
\end{lemma}

\textit{Proof.} We shall use shorthand notation of Lemma 3: $f_t
:= f(x_t)$ and $\varphi_t := \varphi(x_t)$. According to
(\ref{eq21}), one has
 $$
f_t - f_{t-1} \le -\gamma_{t-1} \varphi_{t-1} (\varphi_{t-1} +
\xi_t) + {M\over 2}\, \gamma_{t-1}^2(\varphi_{t-1} + \xi_t)^2 \le
  $$
 $$
\le -\gamma_{t-1} \varphi_{t-1} \xi_t + M \gamma_{t-1}^2
(\varphi_{t-1}^2 + \xi_t^2).
 $$
This implies that $f_t - f_1 \le \Summa_t' + \Summa_t''$, with
 $$
\Summa_t' = \big| \sum_{i=2}^{t} \gamma_{i-1} \varphi_{i-1} \xi_i
\big|, \ \ \ \ \ \Summa_t'' = M \sum_{i=2}^{t} \, \gamma_{i-1}^2
(\varphi^{2}_{i-1} + \xi_i^2).
 $$
Using Lemma 8, one gets
$$
\P (\tau_v < \infty)\, \le\, p(v/w) + P' + P'',
$$
where
$$
P' = \P (\Summa'_{\tau_v} \ge {\delta}/{2}) \ \ \text{ and } \ \
P'' = \P (\Summa''_{\tau_v} \ge {\delta}/{2}).
$$
According to the Chebyshev inequality,
$$
P'\, \le\, \frac{4}{\delta^2}\, E\Summa'^2_{\tau_v} =
\frac{4}{\delta^2} \sum_{i,j=1}^{\infty} E_{ij},
 $$
where
 $$
E_{ij}\, =\, E \left[ \gamma_{i-1} \varphi_{i-1} \xi_i\, \I(i-1 <
\tau_v) \cdot \gamma_{j-1} \varphi_{j-1} \xi_j\, \I(j-1 < \tau_v)
\right].
$$
Using that the values $\gamma_i$,\, $\varphi_i$,\, $\xi_i$, and
$\I (i < \tau_v)$ are $\mathcal{F}_{i}$-measurable, and using
assumptions \A{A1} and \A{A2}, one obtains that for $i \ne j$,\,
$E_{ij} = 0$, and for $i = j$,
 $$
E_{ii} = E \left[ \gamma_{i-1}^2 \varphi_{i-1}^2 \I(i-1 < \tau_v)
\cdot \xi_i^2 \right] \le v^2 e^{-2bi} \sup_{x\in\O_0}
\varphi^2(x) \cdot S.
 $$
Therefore,
 $$
P'\, \le\, \frac{4}{\delta^2} \sum_{i=2}^\infty E_{ii} \le
\frac{4v^2 S}{\delta^2}\ \frac{e^{-4b}}{1 - e^{-2b}}\
\sup_{x\in\O_0} \varphi^2(x).
 $$

Similarly,
$$
P''\, \le\, \frac{2}{\delta}\, E\Summa_{\tau_v}'' =
\frac{2M}{\delta} \sum_{i=2}^\infty E \left[ \gamma^2_{i-1}
(\varphi^{2}_{i-1} + \xi_i^2) \I(i-1 < \tau_v) \right] \le
$$
 $$
\le \frac{2M v^2}{\delta} \sum_{i=2}^\infty e^{-2bi} \left(
\sup_{x\in\O_0} \varphi^2(x) + S \right) = \frac{2M v^2}{\delta}\
\frac{e^{-4b}}{1 - e^{-2b}} \left( \sup_{x\in\O_0} \varphi^2(x) +
S \right).
 $$
Taking
 $$
K\, =\, \left[ \frac{4S}{\delta^2} \sup_{x\in\O_0} \varphi^2(x)\,
+\, \frac{2M}{\delta} \left( \sup_{x\in\O_0} \varphi^2(x) + S
\right) \right] \frac{e^{-4b}}{1 - e^{-2b}},
 $$
one gets that $P' + P'' \le K\, v^2$. Lemma 9 is proved.
 \hfill$\Box$

\begin{lemma}
If $\coefficient > k_+(0)$ then $\sum_t \gamma_t < \infty$.
\end{lemma}

\textit{Proof.} From the definition of $\tau_v$ one easily sees
that if $\tau_v = \infty$ for some $v > 0$, then $\sum_t \gamma_t
< \infty$. This implies that for any $v > 0$
 \begin{equation}\label{1point}
\P \left(\sum \gamma_t\, =\, \infty \right) \le \P(\tau_v =
\infty).
 \end{equation}
Further, by virtue of Lemma 9, if $x_0 \in \O'$ and $\gamma_0 < w$
then
 \begin{equation}\label{2points}
\P (\tau_{\sqrt{w}}\, <\, \infty)\, \le\, Kw + p(1/\sqrt{w}).
 \end{equation}
Combining (\ref{1point}) and (\ref{2points}), one gets that for
any $w > 0$
 \begin{equation}\label{3points}
\P \left(\sum \gamma_t = \infty\ |\ x_0 \in \O' \text{ and }
\gamma_0 < w\right)\, \le\, Kw\, +\, p(1/\sqrt{w}).
 \end{equation}
Define the event $\AA_w = \{ \text{ for some } t,\ x_t \in \O'
\text{ and } \gamma_t < w \}$, then by virtue of (\ref{3points}),
 \begin{equation}\label{101}
\P \left(\sum \gamma_t = \infty\ \big|\ \AA_w \right) \le Kw +
p(1/\sqrt{w}).
 \end{equation}
Denote by $\bar\AA_w$ the complementary event, $\bar\AA_w = \{
\text{ for any } t,\ x_t \not\in \O' \text{ or } \gamma_t \ge w
\}$. By virtue of Lemma 7,
 \begin{equation}\label{102}
\P \left(\sum \gamma_t = \infty \ \, \& \, \ \bar\AA_w \right) =
0.
 \end{equation}
Using (\ref{101}) and (\ref{102}), one gets
 \begin{eqnarray*}
\P \left( \sum \gamma_t =  \infty \right) = \P \left( \sum
\gamma_t = \infty \, \ \& \ \, \AA_w \right) + \P \left( \sum
\gamma_t = \infty \ \, \& \ \, \bar\AA_w \right) \le
 \end{eqnarray*}
 $$
\le  (Kw + p(1/\sqrt{w})) \cdot \P (\AA_w).
 $$
Taking into account that $w$ can be chosen arbitrarily small and
that $Kw + p(1/\sqrt{w}) \to 0$ as $w \to 0^+$, one concludes that
$\P \left( \sum_t \gamma_t =  \infty \right) = 0$.
 \hfill$\Box$
 \vspace{2mm}

Now, we are in a position to prove the theorem. Suppose that
$\coefficient < \inf_z k_-(z)$, then $V_-^{[\coefficient]} =
\emptyset$, and by Lemma 2, $\{ x_t \}$ diverges. So, the
statement (b) of Theorem is proved.

On the other hand, according to Lemma 10, if $\coefficient >
k_+(0)$ then $\sum_t \gamma_t < \infty$, and by Lemmas 1 and 2,
the sequence $\{ x_t \}$ converges to a point from
$V_-^{[\coefficient]}$. Thus, the statement (a) of theorem is also
established.


\section*{Acknowledgements}

This work was partially supported by the R\&D Unit CEOC (Center
for Research in Optimization and Control). The second author (PC)
also gratefully acknowledges the financial support by the
Portuguese program PRODEP   `Medida 5 - Ac\c c\~ao 5.3 - Forma\c
c\~ao Avan\c cada de Docentes do Ensino Superior - Concurso nr.
2/5.3/PRODEP/2001'.


\end{document}